\documentclass[11pt,a4paper]{article}
\usepackage{tikz}
\usetikzlibrary{patterns}

\usepackage{epsfig}
\usepackage{graphics}
\usepackage{amsmath,amssymb}
\usepackage{subfigure}
\usepackage{graphicx}
\usepackage{epstopdf}
\usepackage{hyperref}
\usepackage{pdfsync}
\usepackage{color}

\newtheorem{thm}{Theorem}

\newtheorem{proposition}{Proposition}
\newtheorem{definition}{Definition}

\numberwithin{equation}{section}
\numberwithin{thm}{section}
\numberwithin{lemma}{section}
\numberwithin{proposition}{section}
\numberwithin{definition}{section}
\numberwithin{remark}{section}
\numberwithin{figure}{section}
\numberwithin{table}{section}

\def\E{\mathbb{E}}
\def\Pr{\mathbb{P}}

\title{Minimax signal detection under weak noise assumptions\footnote{This work was supported by the LABEX MILYON (ANR-10-LABX-0070) of Universit\'{e} de Lyon, within the program `Investissements d'Avenir" (ANR-11-IDEX- 0007) operated by the French National Research Agency (ANR).} }
\author{{\em Cl\'ement ~Marteau}, \\
         Univ Lyon, Universit\'e Claude Bernard Lyon 1,\\     CNRS UMR 5208, Institut Camille Jordan, \\
         43 blvd. du 11 novembre 1918, \\F-69622 Villeurbanne cedex, France\\
        {\em Email}:~\texttt{marteau@math.univ-lyon1.fr}
          \\ \\
        and
          \\ \\
        {\em Theofanis ~Sapatinas},\\
        Department of Mathematics and Statistics,\\
        University of Cyprus,\\
        P.O. Box 20537,
        CY 1678 Nicosia,
        Cyprus.\\
        {\em Email}:~\texttt{fanis@ucy.ac.cy}}

\begin{document}
\maketitle

\begin{abstract}
We consider minimax signal detection in the sequence model. Working with certain ellipsoids in the space of square-summable sequences of real numbers, with a ball of positive radius removed,  we obtain upper and lower bounds for the minimax separation radius in the non-asymptotic framework, i.e., for a fixed value of the involved noise level. We use very weak assumptions on the noise (i.e., fourth moments are assumed to be uniformly bounded). In particular, we do not use any kind of Gaussianity or independence assumption on the noise.  It is shown that the established minimax separation rates are not faster than the ones obtained in the classical sequence model (i.e., independent standard Gaussian noise) but, surprisingly, are of the same order as the minimax estimation rates in the classical setting. Under an additional condition on the noise, the classical minimax separation rates are also retrieved in benchmark well-posed and ill-posed inverse problems. 


\end{abstract}

\newpage

\section{Introduction}
We consider the following sequence model (SM),
\begin{equation}\label{1.0.0}
y_k=b_k \theta_k+\varepsilon \, \xi_k, \quad k \in \mathcal{N},
\end{equation}
where $\mathcal{N}$ can be either $\mathbb{N}=\{1,2,\ldots\}$ or  $\mathbb{N}_n=\lbrace1,\dots, n \rbrace$ for some $n \geq 1$,  $b=(b_k)_{k \in \mathcal{N}}$ is a known positive sequence,  $\theta=(\theta_k)_{k\in\mathcal{N}} \in l^2(\mathcal{N})$ is the unknown signal of interest, $\xi=(\xi_k)_{k \in\mathcal{N}}$ is a sequence of random variables (the noise), and $\varepsilon >0$ is a known parameter (the noise level). The observations are given by the sequence $y=\{y_k\}_{k \in \mathcal{N}}$ from the SM (\ref{1.0.0}) and their joint law is denoted by $\Pr_{\theta,\xi}$. Here, $l^2(\mathcal{N})$ denotes the space of squared-summable sequence of real numbers, i.e., 
$$
l^2(\mathcal{N}) = \left\lbrace \theta \in \mathbb{R}^\mathcal{N}:\; \|\theta\|^2:=\sum_{k \in \mathcal{N}}\theta_j^2 < +\infty \right\rbrace.
$$
Let $\mathcal{C}>0$ be a known fixed constant.  Concerning the noise, we will assume that $\xi \in \Xi$, where 
\begin{equation}
\Xi := \Xi(\mathcal{C})=  \left\lbrace \xi: \ \E[\xi_k] = 0, \ \E[\xi_k^2]= 1 \ \forall k \in \mathcal{N} \ \mathrm{and} \ \sup_{k \in \mathcal{N}} \E[\xi_k^4] \leq \mathcal{C} < + \infty \right\rbrace.
\label{eq:Xi}
\end{equation}

The SM (\ref{1.0.0}) arises in many well-known situations. Consider for instance the stochastic differential equation
$$ dZ_\varepsilon(t) = A f(t) + \epsilon d U(t), \quad t\in [0,1],$$
where $A$ is a known bounded linear operator acting on $L^2([0,1])$, $f(\cdot) \in L^2([0,1])$  is the unknown response function that one wants to detect or estimate, $U(\cdot)$ is a given stochastic process on $[0,1]$ and $\varepsilon>0$ is a known parameter (the noise level).  For the sake of simplicity, we only consider the case when $A$ is injective (meaning that $A$ has a trivial nullspace).

\begin{itemize}
\item[$\bullet$] Let $U(\cdot)=W(\cdot)$ be the standard Wiener process. Then, if $A$ is the identity operator, we can retrieve the SM (\ref{1.0.0}) in the Fourier domain with $b_k=1$ for all $k \in \mathcal{N} = \mathbb{N}$ and the $\xi_k$, $k \in \mathcal{N}$, are independent standard Gaussian random variables (direct problem).  If $A$ is a self-adjoint operator with an eigen-decomposition, we can retrieve the SM (\ref{1.0.0}) where $b_k > b_0$ for some $b_0>0$ for all $k \in \mathcal{N} = \mathbb{N}$ and the $\xi_k$, $k \in \mathcal{N}$,  are independent standard Gaussian random variables (well-posed inverse problems). If $A$ is a compact operator, we can retrieve the SM (\ref{1.0.0}) where $b_k>0$ for all $k \in\mathcal{N}$ (since $A$ is injective) with $b_k \rightarrow 0$ as $k \rightarrow + \infty$ and the $\xi_k$, $k \in \mathcal{N}$, are independent standard Gaussian random variables (ill-posed inverse problems). For more details regarding all these models, we refer to, e.g., \cite{Cavalier_book}.\\
\item[$\bullet$] Let $U(\cdot)= W_{-\gamma}(\cdot)$, $\gamma \in ]0,1/2[$, be the truncated fractional Brownian motion and let $A$ be the identity operator. Then, we can retrieve the SM (\ref{1.0.0}) in the spline domain with $b_k = (\pi k)^{-2\gamma}(1+o(1))$ as $k\rightarrow + \infty$ and the $\xi_k$, $k \in \mathcal{N}$, are (non-independent) standard Gaussian random variables. For more details, we refer to, e.g., \cite{J_1999}, \cite{Cavalier_frac}.
\end{itemize}

The non-parametric inverse regression problem also provides observations of the form (\ref{1.0.0}). Indeed, consider the model
$$ Z_i = Af\left(  \frac{i}{n} \right) + \frac{1}{\sqrt{n}}\, \eta_i, \quad i \in \lbrace 1, \dots, n \rbrace,$$
where $A$ is a known (injective) bounded linear operator acting on $L^2([0,1])$, $f(\cdot) \in L^2([0,1])$  is the unknown response function that one wants to detect or estimate, and $\eta_i$, $i \in \mathbb{N}_n$, is a sequence of independent and identically distributed random variables with zero mean, variance one and finite fourth moment. Given any appropriate bases (or, even, a tight frame, see, e.g., \cite{Mallat}, p. 126), we can retrieve the SM (\ref{1.0.0}) with $b_k=1$ for all $k \in \mathcal{N}_n$ when $A$ is the identity operator, see, e.g., \cite{Tsybakov}, Chapter 1.  When $A$ is a compact operator, we can retrieve an approximation of the SM (\ref{1.0.0}) where $(b_k)_{k \in \mathbb{N}_n}$ is a fixed sequence that depends on $A$, see, e.g., \cite{Munk_testing}.  \\

Minimax signal detection has been considered in the literature over the last two decades. We refer to, e.g., \cite{Baraud}, \cite{LLM_2012}, \cite{IS_2003}, \cite{LLM_2011}, \cite{ISS_2012}, \cite{ISS_2011}, \cite{MM_2013}. All these contributions consider the classical Gaussian sequence model (\ref{1.0.0}), i.e,. where the $\xi_k$, $k \in \mathcal{N}$, are independent standard Gaussian random variables. We refer to \cite{MS_2014} for a survey on available results and a discussion on the link between asymptotic (the noise level is assumed to tend to zero) and non-asymptotic  (the noise level is assumed to be fixed) approaches to minimax signal detection. The aim of this work is to obtain upper and lower bounds on the minimax separation radius in the non-asymptotic framework, for the general model (\ref{1.0.0}) under weak assumptions on the noise, i.e., when $\xi \in \Xi$, where the set $\Xi$ has been introduced in (\ref{eq:Xi}).  In particular, we do not use any kind of Gaussianity or independence assumption on the noise. We prove that the minimax separation rates are not faster than the ones obtained in the classical sequence model (see, e.g., \cite{Baraud}, \cite{LLM_2012}, \cite{IS_2003}, \cite{LLM_2011}, \cite{ISS_2012}, \cite{ISS_2011}, \cite{MM_2013}) but,  surprisingly, are of the same order as the minimax estimation rates in the classical setting. Moreover, under additional conditions on the noise, we show that the classical minimax separation rates can be retrieved in benchmark well-posed and ill-posed inverse problems.      \\

Throughout the paper, we use the following notations.  Given two sequences $(c_k)_{k \in \mathcal{N}}$ and $(d_k)_{k \in \mathcal{N}}$ of real numbers, $c_k \sim d_k$ means that there exist $0<\kappa_0 \leq \kappa_1 <\infty$ such that $\kappa_0 \leq c_k/d_k \leq \kappa_1$ for all $k \in \mathcal{N}$, while $c_k\lesssim d_k$ (resp. $c_k \gtrsim d_k$) means $c_k \leq c_0 \; d_k$ (resp. $c_k \geq c_0 \; d_k$) for some $c_0>0$ for all $k \in \mathcal{N}$. Also, $ x\wedge y:=\min(x,y)$, for all $x,y \in \mathbb{R}$.

\section{Minimax Signal Detection}
\label{msd}
Given observations from the SM (\ref{1.0.0}), we consider the signal detection problem, i.e., our aim is to test 
\begin{equation}
H_0: \theta=0 \;\; \mathrm{versus} \;\; \ H_1: \theta  \in \Theta_a(r_\varepsilon).
\label{testing_pb0}
\end{equation}
Given a non-decreasing sequence $a=(a_k)_{k\in \mathcal{N}}$ of positive real numbers, with $a_k \rightarrow +\infty$ as $k \rightarrow +\infty$ when $\mathcal{N}=\mathbb{N}$, and a radius $r_\varepsilon>0$, the set $\Theta_a(r_\varepsilon)$ is defined as
\begin{equation}
\label{def:f-set}
\Theta_a(r_\varepsilon) = \left\lbrace \theta \in \mathcal{E}_a, \ \|\theta \| \geq r_\varepsilon \right\rbrace,
\end{equation}
where
$$ 
\mathcal{E}_{a} = \left\lbrace \theta\in l^2(\mathcal{N}), \ \sum_{k\in \mathcal{N}} a_k^2 \theta_k^2 \leq  1\right\rbrace.
$$
The set $ \mathcal{E}_{a}$ can be seen as a condition on the decay of $\theta$. The cases where the sequence $a$ increases very fast correspond to the signal $\theta$ with small coefficients. In such a case, the corresponding signal can be considered as being `smooth'. The sequence $a$ being fixed, the main issue for the minimax signal detection problem (\ref{testing_pb0})-(\ref{def:f-set}) is then to characterize the values of the radius $r_\varepsilon >0$ for which both hypotheses $H_0$ (called the null hypothesis) and $H_1$ (called the alternative hypothesis) are `separable'.\\

In the following, a (non-randomized) test $\Psi:=\Psi(y)$ will be defined as a measurable function of the observation $y=(y_k)_{k\in\mathcal{N}}$ from the SM (\ref{1.0.0}) having values in the set $\lbrace 0,1 \rbrace$. By convention, $H_0$ is rejected if $\Psi=1$ and  $H_0$ is not rejected if $\Psi=0$. Then, given a test $\Psi$, we can investigate 
\begin{itemize}
\item the type I (first kind) error probability defined as 
\begin{equation}
\sup_{\xi \in \Xi} \mathbb{P}_{0,\xi}( \Psi =1),
\label{eq:type1-F}
\end{equation} 
which measures the worst probability of rejecting $H_0$ when $H_0$ is true (i.e., $\theta=0$, $\xi \in \Xi$, where $\Xi$ is defined in (\ref{eq:Xi})); it is often constrained as being bounded by a prescribed level $\alpha \in ]0,1[$, and
\item the type II (second kind) error probability defined as 
\begin{equation}
\sup_{\substack{\theta\in \Theta_a(r_\varepsilon) \\ \xi \in \Xi}} \Pr_{\theta,\xi}(\Psi=0),
\label{eq:type2}
\end{equation}
which measures the worst possible probability of not rejecting $H_0$ when $H_0$ is not true (i.e., when $\theta \in \Theta_a(r_\varepsilon)$ and $\xi \in \Xi$, where $\Xi$ is defined in (\ref{eq:Xi})); one would like to ensure that it is bounded by a prescribed level $\beta \in ]0,1[$.
\end{itemize}
\vspace{0.2cm}

We emphasize that in the classical minimax signal detection problem, the protection against all possible noise distributions (i.e., $\sup_{\xi \in \Xi}$) is not required, since the noise distribution is completely known. However, in the more general setting that we consider, in order to produce some kind of robustness, we have adapted the definitions of type I and type II error probabilities to accommodate the (possible) uncertainty on the noise.  \\

Let $\alpha,\beta\in ]0,1[$ be given, and let $\Psi_\alpha$ be an $\alpha$-level test, i.e., $\Psi_\alpha$ is such that $\sup_{\xi\in \Xi}\mathbb{P}_{0,\xi}( \Psi_\alpha =1) \leq \alpha$.

\begin{definition}
The separation radius of the $\alpha$-level test $\Psi_\alpha$ over the class $\mathcal{E}_a$ is defined as
$$ r_\varepsilon(\mathcal{E}_a,\Psi_\alpha,\beta) := \inf \left\lbrace r_\varepsilon>0: \ \sup_{\substack{\theta\in \Theta_a(r_\varepsilon) \\ \xi \in \Xi}} \Pr_{\theta,\xi}(\Psi_\alpha=0)  \leq \beta\right\rbrace.$$
\label{eq:minimax_separation_radius}
\end{definition}

\vspace{-0.6cm}

In some sense, the separation radius $r_\varepsilon(\mathcal{E}_a,\Psi_\alpha,\beta)$ corresponds to the smallest possible value of the available signal $\| \theta \|$ for which $H_0$ and $H_1$ can be `separated' by the $\alpha$-level test $\Psi_\alpha$ with prescribed type I and type II error probabilities, $\alpha$ and $\beta$, respectively.

\begin{definition}
\label{def:minimax_radius}
The minimax separation radius $\tilde{r}_\varepsilon:=\tilde{r}_{\varepsilon}(\mathcal{E}_a, \alpha, \beta)>0$ over the class $\mathcal{E}_a$ is defined as
\begin{equation}
\tilde r_{\varepsilon}:= \inf_{\tilde \Psi_\alpha} r_\varepsilon(\mathcal{E}_a, \tilde\Psi_\alpha,\beta),
\label{eq:minimax_radius}
\end{equation}
where the infimum is taken over all $\alpha$-level tests $\tilde \Psi_\alpha$.
\end{definition}

The minimax separation radius $\tilde r_{\varepsilon}$ corresponds to the smallest radius $r_{\varepsilon} >0$ such that there exists some $\alpha$-level test $\tilde \Psi_\alpha$ for which the type II error probability is not greater than $\beta$.\\

It is worth mentioning that Definitions \ref{eq:minimax_separation_radius} and \ref{def:minimax_radius} are valid for any {\em fixed} $\varepsilon >0$ (i.e., it is not required that $\varepsilon \rightarrow 0$). The performances of any given test $\Psi_\alpha$ is easy to handle in the sense that the type I error probability is bounded by $\alpha$ (i.e., $\Psi_\alpha$ is an $\alpha$-level test), and that the dependence of the minimax separation radius $\tilde r_{\varepsilon}$ with respect to given $\alpha$ and $\beta$ can be precisely described.

\section{Control of the Upper and Lower bounds}
\subsection{The spectral cut-off test and control of the upper bound}

We define below a spectral cut-off test for the SM model (\ref{1.0.0}) with $\xi \in \Xi$, where $\Xi$ is defined in (\ref{eq:Xi}). First, we show that it is an $\alpha$-level test and then we obtain an upper bound for its type II error probability.\\

Given a bandwidth $D \in \mathcal{N}$ and $\alpha \in ]0,1[$, we consider the following spectral cut-off test
\begin{equation}
\Psi_{\alpha,D} := \mathbf{1}_{\lbrace  T_D \geq t_{1-\alpha,D}  \rbrace},
\label{eq:SCOT-fanis}
\end{equation}
where 
$$ 
T_D = \sum_{k=1}^D b_k^{-2} (y_k^2 - \varepsilon^2)
$$
and $t_{1-\alpha,D}$ denotes a threshold depending on $\alpha$ and $D$. It is easily seen that, for all $D \in \mathcal{N}$,
$$ \E_{\theta, \xi} [ T_D] = \sum_{k=1}^D \theta_k^2,$$
and 
\begin{equation} 
\mathrm{Var}_{0,\xi}(T_D)  = \underbrace{\varepsilon^4 \sum_{k=1}^D b_k^{-4}  \E [ (\xi_k^2 -1)^2]}_{:= R_0(D)} +  \underbrace{ \varepsilon^4 \sum_{\substack{k,l=1 \\  k \not = l}}^D b_k^{-2} b_l^{-2}  \E [ (\xi_k^2 -1)  (\xi_l^2 -1)]}_{:=S_0(D)},
\label{eq:Var2}
\end{equation}
where the assumption $\xi \in \Xi$ guarantees that the above variance is finite for every $D \in \mathcal{N}$.

\begin{proposition}
\label{prop:type1}
Let $\alpha \in ]0,1[$ be given. Consider the spectral cut-off test $\Psi_{\alpha,D}$ defined in (\ref{eq:SCOT-fanis}). Then, for all $\varepsilon >0$,
$$ \sup_{\xi \in \Xi} \Pr_{0,\xi} ( \Psi_{\alpha,D} =1) \leq \alpha$$
as soon as 
\begin{equation}
t_{1-\alpha,D} \geq \frac{1}{\sqrt{\alpha}}  \sqrt{ R_0(D) + S_0(D)}.
\label{eq:cond_t}
\end{equation}
\end{proposition}

\noindent
The proof of this proposition is postponed to Section \ref{s:ptype1}.\\

\noindent
\textbf{Remarks:} 
\begin{itemize}
\item Using simple bounds, it is easily seen that
\begin{equation}
R_0(D) + S_0(D) \leq  C_{1} \varepsilon^4  \sum_{1\leq k \leq D} b_k^{-4} + C_{1} \varepsilon^4 \left( \sum_{1\leq k \leq D} b_k^{-2} \right)^2  \leq 2C_{1} \varepsilon^4 \left( \sum_{1\leq k \leq D} b_k^{-2} \right)^2 , 
\label{eq:range}
\end{equation}
where
\begin{equation}
C_{1} := \sup_{\xi\in \Xi} \ \sup_{k\in \mathcal{N}} \E [ (\xi_k^2-1)^2]< +\infty,
\label{eq:C_one}
\end{equation}
since
$$ \sum_{k=1}^D b_k^{-4} \leq (\max_{1\leq k \leq D} b_k^{-2} )  \sum_{k=1}^D b_k^{-2} \leq \left( \sum_{k=1}^D b_k^{-2} \right)^2.$$
Hence, the choice
\begin{equation}
t_{1-\alpha,D} = K_1 \varepsilon^2 \sum_{k=1}^D b_k^{-2}, \quad \mathrm{where} \quad K_1 = \frac{\sqrt{2C_{1}}}{\sqrt{\alpha}},
\label{eq:t2}
\end{equation}
ensures that (\ref{eq:cond_t}) is satisfied and that the spectral cut-off test $\Psi_{\alpha,D}$ defined in (\ref{eq:SCOT-fanis}) is an $\alpha$-level test.  
\item In the classical setting (i.e., independent Gaussian noise), the threshold $t_{1-\alpha,D}$ can be chosen as the $(1-\alpha)$-quantile of the variable $T_D$ under $H_0$. This is no more the case here since only a uniform bound on the fourth moment of the sequence $\xi_k$, $k \in \mathcal{N}$, is available. 
\end{itemize}

\begin{proposition}
\label{prop:upperbound}
Let $\alpha,\beta \in ]0,1[$ be given. Consider the spectral cut-off test $\Psi_{\alpha,D}$ defined in (\ref{eq:SCOT-fanis}).  Select the threshold $t_{1-\alpha,D}$ as in  (\ref{eq:t2}). Then, for all $\varepsilon >0$,
$$ \sup_{\substack{\theta\in \Theta_a(r_{\varepsilon,D}) \\ \xi \in \Xi}} \Pr_{\theta,\xi}(\Psi_{\alpha,D}=0) \leq \beta,$$
for all radius $r_{\varepsilon,D}>0$ such that 
$$ r_{\varepsilon,D} \geq \mathcal{C}_\beta \varepsilon^2 \sum_{k=1}^D b_k^{-2} + a_D^{-2},$$
where $\mathcal{C}_\beta>0$ is the solution of the equation (\ref{eq:condCb}).
\end{proposition}

\noindent
The proof of Proposition is postponed to Section \ref{s:pupperbound}. \\

\noindent
\textbf{Remark:} For practical purposes, the solution $\mathcal{C}_\beta>0$ of equation (\ref{eq:condCb}) can be chosen as $ \mathcal{C}_\beta = 8 K_2/ \beta$. In particular, there exists some $\beta_0> 0$ such that $1- K_1 \mathcal{C}_\beta^{-1} \geq 1/2$ for all $\beta \leq \beta_0$, hence ensuring that (\ref{eq:Cb}) is satisfied for all $\beta$ small enough. 

\subsection{Control of the lower bound}

We propose below a lower bound on the minimax type II error probability for the SM (\ref{1.0.0}) with $\xi \in \Xi$, where $\Xi$ is defined in (\ref{eq:Xi}). In the sequel, the term $\inf_{\Psi_\alpha}$ corresponds  to an infimum taken over all possible $\alpha$-level tests.

\begin{proposition}
\label{prop:lowerbound}
Let $\alpha \in ]0,1[$ and $\beta \in ]0,1-\alpha[$ be fixed. Then, for all $\varepsilon >0$
$$ \inf_{\Psi_\alpha} \sup_{\substack{\theta \in \Theta_a(r_{\varepsilon}) \\ \xi \in \Xi}} \Pr_{ \theta,\xi}(\Psi_\alpha =0) \geq \beta,$$
for all $D \in \mathcal{N}$ and $r_{\varepsilon}>0$ such that 
$$ r_\varepsilon^2 \leq \left( \frac{1}{4} \ln( \mathcal{C}_{\alpha,\beta}) \right)\varepsilon^2 \sum_{k=1}^D b_k^{-2} \wedge a_D^{-2},$$
where $\mathcal{C}_{\alpha,\beta}=1 + 4 (1-\alpha - \beta)^2$.
\end{proposition}

\noindent
The proof of Proposition \ref{prop:lowerbound} is postponed to Section \ref{s:plowerbound}. The main difficulty is to construct an appropriate distribution for $\xi$ that will allow one to obtain the largest possible lower bound. 


\section{Minimax Separation Radius}

The following theorem provides upper and lower bounds for the minimax separation radius $\tilde r_\varepsilon>0$ in the SM (\ref{1.0.0}) with $\xi \in \Xi$, where $\Xi$ is defined in (\ref{eq:Xi}).

\begin{thm}
\label{thm:1}
Let $\alpha, \beta \in\, ]0,1[$ be given. Then, for all $\varepsilon >0$, the minimax separation radius $\tilde r_\varepsilon>0$ satisfies 
\begin{equation}
\sup_{D\in \mathcal{N}} \left[ \left( \frac{1}{4} \ln( \mathcal{C}_{\alpha,\beta}) \right) \varepsilon^2 \sum_{k=1}^D b_k^{-2} \wedge a_D^{-2} \right]  \leq \tilde r_\varepsilon^2 \leq \inf_{D\in\mathcal{N}} \left[ \mathcal{C}_\beta   \varepsilon^2 \sum_{k=1}^D b_k^{-2} + a_D^{-2}\right] ,
\label{eq:main}
\end{equation}
where $\mathcal{C}_\beta>0$ is the solution of the equation (\ref{eq:condCb}) and $\mathcal{C}_{\alpha,\beta}=1 + 4 (1-\alpha - \beta)^2$. 
\end{thm}

\noindent
The proof of Theorem \ref{thm:1} is postponed to Section \ref{s:pthm1}.\\

\noindent
\textbf{Remark:} If both sequences $a=(a_k)_{k\in \mathcal{N}}$ and $b^{-1}=(b_k^{-1})_{k\in\mathcal{N}}$ are non-decreasing and satisfy
\begin{equation}
a_\star \leq \frac{a_{D-1}}{a_{D}} \leq a^\star \quad \mathrm{and} \quad b_\star \leq \frac{b_{D-1}}{b_{D}} \leq b^\star \quad \mbox{for all} \quad D\in \mathcal{N} \setminus \lbrace 1 \rbrace,
\label{eq:hyp_ab}
\end{equation}
for some constants $0<a_\star \leq a^\star <\infty$ and $0<b_\star \leq b^\star <\infty$, then, it is easily seen that both upper and lower bounds on the minimax separation radius $\tilde r_\varepsilon>0$, established in Theorem \ref{thm:1}, are of the same order. This follows easily working along the same lines of the proof of Proposition 4.1 in \cite{MS_2014}. We note also that the condition (\ref{eq:hyp_ab}) is satisfied for various combinations of interest, among them: (i) mildly ill-posed inverse problems ($b_k \asymp k^{-t}$, $k\in\mathbb{N}$, for some $t>0$) with ordinary smooth functions ($a_k  \asymp k^{s}$, $k\in\mathbb{N}$, for some $s>0$), (ii) severely ill-posed inverse problems ($b_k \asymp e^{-kt}$, $k\in\mathbb{N}$, for some $t>0$) with ordinary smooth functions ($a_k \asymp k^{s}$, $k\in\mathbb{N}$, for some $s>0$), and (iii) mildly ill-posed inverse problems ($b_k \asymp k^{-t}$, $k\in\mathbb{N}$, for some $t>0$) with super-smooth functions ($a_k \asymp e^{ks}$, $k\in\mathbb{N}$, for some $s>0$). Among the possible situations where the condition (\ref{eq:hyp_ab}) is not satisfied, one can mention, for instance, power-exponential behaviors ($a_k \asymp e^{k^l s}$, $j\in\mathbb{N}$, for some $s>0$ and $l>1$, or $b_k \asymp e^{-k^rt}$, $k\in\mathbb{N}$, for some $t>0$ and $r>1$). See also Remark 4.3 in \cite{MS_2014}.\\

\noindent
\textbf{Remark:} Note that the upper and lower bounds on the minimax separation radius $\tilde r_\varepsilon>0$, established in Theorem \ref{thm:1}, are quite different compared to the classical minimax separation radii available in the literature, obtained in the SM (\ref{1.0.0}) with independent standard Gaussian noise (see, e.g., \cite{MS_2014}).  Although the bias terms $a_D^{-2}$ coincide, the corresponding variance terms differ. In particular, in the SM (\ref{1.0.0}) with $\xi \in \Xi$, where $\Xi$ is defined in (\ref{eq:Xi}), the variance term is of order $\varepsilon^2 \sum_{k=1}^D b_k^{-2}$, while for the SM (\ref{1.0.0}) with independent standard Gausiian noise, the variance term is of order $\varepsilon^2 \sqrt{\sum_{k=1}^D b_k^{-4}}$.  We stress that the term $\varepsilon^2 \sum_{k=1}^D b_k^{-2}$ is not greater than the term$\varepsilon^2 \sqrt{\sum_{k=1}^D b_k^{-4}}$, which entails that the minimax separation rates are not faster compared to the ones obtained in the classical model. It is also worth mentioning that, surprisingly,  the bias and variance terms in the SM (\ref{1.0.0}) with $\xi \in \Xi$, where $\Xi$ is defined in (\ref{eq:Xi}), are of the same order of the corresponding terms in the classical minimax estimation setting. In particular, the minimax separation rates in our general setting coincide with the minimax estimation rates obtained in the classical estimation setting. For illustrative purposes, the Table \ref{tab:sep} (see also Table 1 in \cite{Cavalier_book}) provides these minimax separation rates for benchmark problems, i.e., well-posed, mildly ill-posed and severely ill-posed problems for ellipsoids with ordinary smooth and super-smooth sequences.  \\

\begin{table}[h]
\centering
\begin{tabular}{|c|c|c|}
\hline
{ \textbf{Minimax separation}  }&ordinary-smooth &super-smooth \\
{ \textbf{rate} $(\tilde r_\varepsilon^2)$  } & $a_k \sim k^{s}$ & $a_k\sim \exp\{ks\}$ \\
\hline 
well-posed  & $\varepsilon^{4s/(2s+1)}$ & $\varepsilon^2(\ln \varepsilon^{-1})$\\
$ b_k \sim 1$ & & \\
\hline
mildly ill-posed &$\varepsilon^{4s/(2s+2t+1)}$&$
\varepsilon^2(\ln \varepsilon^{-1})^{2t+1}$\\
 $b_k\sim k^{-t}$ & & \\
\hline
severely ill-posed & $(\ln \varepsilon^{-1})^{-2s}$&$
\varepsilon^{4s/(2s+2t)}$\\ 
$b_k \sim \exp\{-kt\}$ & & \\
\hline
\end{tabular}
\caption{\textit{Minimax separation rates for the SM (\ref{1.0.0}) with $\xi \in \Xi$, where $\Xi$ is defined in (\ref{eq:Xi}).}}
\label{tab:sep}
\end{table}

\noindent
\textbf{Remark:} If the supremum over all possible noise distributions $\xi \in \Xi$ is not considered in the definition of  type I and type II error probabilities, then it is easily seen that upper bound on the type II error probability obtained in Proposition \ref{prop:upperbound} still holds true. However, the corresponding lower bound obtained in Proposition \ref{prop:lowerbound} is only true under Gaussianity. This implies that the minimax separation rates displayed in Table \ref{tab:sep} are still valid in the SM (\ref{1.0.0}) with non-independent standard Gaussian noise $\xi$.

\section{An Additional Condition on the Noise to Obtain the Classical Minimax Separations Rates}

In this section, it is demonstrated that, under an additional condition on the noise $\xi \in \Xi$ in the SM (\ref{1.0.0}), one is able to retrieve the classical minimax separation rates in benchmark well-posed and ill-posed inverse problems. \\

Recall from equation (\ref{eq:var}), displayed in the proof of Proposition \ref{prop:upperbound}, that  the variance of $T_D$ can be written as
$$ \mathrm{Var}_{\theta,\xi}(T_D) = R_\theta(D) + S_\theta(D),$$
where 
$$ R_\theta(D):=  \sum_{k=1}^D b_k^{-4}  \mathrm{Var}_{\theta,\xi} (y_k^2 - \varepsilon^2)$$
and 
\begin{eqnarray*}
S_\theta(D)& := & \sum_{\substack{k,l=1 \\  k \not = l}}^D b_k^{-2} b_l^{-2}  \mathrm{Cov}_{\theta,\xi}(y_k^2- \varepsilon^2,y_l^2-\varepsilon^2).
\end{eqnarray*}
In the classical setting (i.e., independent standard Gaussian noise $\xi$),  $S_\theta(D)=0$ for all $\theta \in l^2(\mathcal{N})$. Hence, in order to retrieve the classical minimax separation rates in the SM (\ref{1.0.0}) with $\xi \in \Xi$, where $\Xi$ is defined in (\ref{eq:Xi}),  $S_\theta(D)$ needs to be of the order of $R_\theta(D)$. We achieve this separately under the null and the alternative hypotheses, for benchmark problems, such as, well-posed, mildly ill-posed and severely ill-posed inverse problems. \\

We stress that in this section, we will only deal with upper bounds. Indeed, the lower bounds established previously in the literature (see, e.g., \cite{MS_2014}, Theorem 4.1) for the independent standard Gaussian noise are still valid in our set-up.

\subsection{Well-posed and mildly ill-posed inverse problems}
We assume that 
$$ b_k \sim k^{- t} \quad \forall k \in \mathcal{N}$$
for some $t\geq 0$ ($t=0$ refers to well-posed inverse problems while $t>0$ refers to mildly ill-posed inverse problems).  We start our discussion under the null hypothesis. Recall from (\ref{eq:Var2}) that 
\begin{eqnarray*}
\mathrm{Var}_{0,\xi}(T_D)  
& = & \varepsilon^4 \sum_{k=1}^D b_k^{-4}  \E [ (\xi_k^2 -1)^2]+  \varepsilon^4 \sum_{\substack{k,l=1 \\  k \not = l}}^D b_k^{-2} b_l^{-2}  \E [ (\xi_k^2 -1) (\xi_l^2 -1)] \\
& = & \varepsilon^4 \sum_{k=1}^D b_k^{-4}  \mathrm{Var}(\xi_k^2)+   \varepsilon^4 \sum_{\substack{k,l=1 \\  k \not = l}}^D b_k^{-2} b_l^{-2}  \mathrm{Cov}(\xi_k^2, \xi_l^2) \\
& := & R_0(D) + S_0(D).
\end{eqnarray*}
Using simple calculations, we can see that 
$$ R_0(D) \sim \varepsilon^4 \sum_{k=1}^D k^{4t}  \sim \varepsilon^4 D^{4t+1}.$$
Our aim is to exhibit a condition for which $S_0(D)$ is (at least) of the same order as $R_0(D)$.\\
\\
\textbf{Assumption $\mathcal{H}_\mathcal{D}$:} \textit{Let $\xi\in \Xi$, where $\Xi$ is defined in (\ref{eq:Xi}),  and, for all $k,l \in \mathcal{N}$, let $(\xi_k,\xi_l)'$ be a bivariate Gaussian random vector. Moreover, there exists $s>0$ such that}
$$ \rho_{kl} := |\mathrm{Cov}(\xi_k,\xi_l)| \lesssim \frac{1}{|k-l|^{s}} \quad \forall k,l \in \mathcal{N}, \ k\not = l.$$

Due to the Isserlis Theorem (see, e.g., \cite{Isserlis}), it can be seen that, thanks to Assumption $\mathcal{H}_\mathcal{D}$, for all $k,l \in \mathcal{N}$, with $k\not = l$, 
\begin{equation}
\mathrm{Cov}(\xi_k^2,\xi_l^2) = 2 \mathrm{Cov}^2(\xi_k,\xi_l) \lesssim \frac{1}{|k-l|^{2s}}
\label{eq:anoise1}
\end{equation}
and 
\begin{equation}
\E [ (\xi_k^2-1) \xi_l] = \E [ (\xi_k^2-1) \xi_l] = 0.
\label{eq:anoise2}
\end{equation}
These results allow us to propose a sharp control of the variance of $T_D$ under the null hypothesis.

\begin{proposition}
\label{prop:short1}
Assume that Assumption $\mathcal{H}_\mathcal{D}$ holds with $s>1/2$. Then, 
$$ S_0(D) = o(R_0(D)) \quad \mathrm{as} \quad D \rightarrow + \infty.$$
\end{proposition}

\noindent
The proof of Proposition \ref{prop:short1} is postponed to Section \ref{s:pshort1}. \\

Now, we propose a similar analysis under the alternative hypothesis.  

\begin{proposition}
\label{prop:short2}
Assume that Assumption $\mathcal{H}_\mathcal{D}$ holds with $s>1/2$. Then,  for all $\gamma \in ]0,1[$,
$$ \mathrm{Var}_{\theta,\xi}(T_D) \lesssim   (1+\gamma^{-1}) \varepsilon^4 \sum_{k=1}^D b_k^{-4}  + \gamma \left( \sum_{k=1}^D \theta_k^2 \right)^2.$$
\end{proposition}

\noindent
The proof of Proposition \ref{prop:short2} is postponed to Section \ref{s:pshort2}. \\

Starting from (\ref{eq:2kind}), and using Propositions \ref{prop:short1} and \ref{prop:short2}, we get 
\begin{eqnarray*}
\Pr_{\theta,\xi} ( \Psi_{\alpha,D} =0)
& \leq & \frac{ \mathrm{Var}_{\theta,\xi}(T_D)}{\left( \sum_{k=1}^D \theta_k^2 - t_{1-\alpha,D}  \right)^2}  \\
& \lesssim & \frac{(1+\gamma^{-1}) \varepsilon^4 \sum_{k=1}^D b_k^{-4}  + \gamma \left( \sum_{k=1}^D \theta_k^2 \right)^2}{\left( \sum_{k=1}^D \theta_k^2 - \varepsilon^4 \sum_{k=1}^D b_k^{-4}  \right)^2} \\
& \leq & \beta
\end{eqnarray*}
provided
$$ \sum_{k=1}^D \theta_k^2 \gtrsim   \varepsilon^2 \sqrt{\sum_{k=1}^D b_k^{-4} },$$
which holds as soon as 
$$ \| \theta \|^2 \gtrsim   a_D^{-2} +  \varepsilon^2 \sqrt{\sum_{k=1}^D b_k^{-4} }.$$
The last inequality provides a classical condition that has been already discussed in, e.g., \cite{Baraud}, \cite{ISS_2012} and \cite{LLM_2012}, or in Theorem 4.1 of \cite{MS_2014}, in the specific case where the noise $\xi$ in the SM (\ref{1.0.0}) is assumed to be independent standard Gaussian. This entails that the Assumption $\mathcal{H}_\mathcal{D}$ suffices to retrieve the classical minimax separation rates for mildly ill-posed inverse problems. 


\subsection{Severely ill-posed inverse problems}

We assume in this section that 
$$ b_k \sim e^{-k t} \quad \forall k \in \mathcal{N}$$
for some $t>0$. Since minimax estimation and minimax separation rates in the classical setting are of the same order (see, e.g., Tables 2 and 3 in \cite{ISS_2012}), we stress that non-independence does not deteriorate the classical minimax separation rates. In other words, the independent standard Gaussian assumption on noise $\xi$ is not needed to get the classical minimax separation rates for severely ill-posed inverse problems.

\section{Concluding Remarks}

We have established minimax separation rates in a general Gaussian sequence model, i.e., the noise need neither to be independent nor standard Gaussian. These rates are not faster than the ones obtained in the classical setting (i.e., independent standard Gaussian noise) but, surprisingly, are of the same order as the minimax estimation rates in the classical setting. The involved spectral cut-off test depends on the unknown smoothness parameter of the signal under the alternative hypothesis. It is therefore of paramount importance in practical applications to provide minimax testing procedures that do not explicitly depend on the associated smoothness parameter. This is, usually, referred to as the `adaptation' problem. However, such an investigation needs careful attention that is beyond the scope of the present work. In particular, the dependency of the involved constant with respect to the level $\alpha$ has a more intricate form than the one involved in the classical setting.

\section{Appendix}

\subsection{Proof of Proposition \ref{prop:type1}}
\label{s:ptype1}

Let $\xi \in \Xi$ be fixed. Using the Markov inequality, we get
$$ \Pr_{0,\xi} ( \Psi_{\alpha,D} =1) = \Pr_{0,\xi}(T_D \geq t_{1-\alpha,D}) \leq \frac{\mathrm{Var}_{0,\xi}(T_D)}{t_{1-\alpha,D}^2} \leq \frac{R_0(D)+ S_0(D)}{t_{1-\alpha,D}^2} \leq \alpha$$
provided 
$$
t_{1-\alpha,D} \geq \frac{1}{\sqrt{\alpha}}  \sqrt{ R_0(D) + S_0(D)}.
$$
\begin{flushright}
$\Box$
\end{flushright}

\subsection{Proof of Proposition \ref{prop:upperbound}}
\label{s:pupperbound}

Let $\xi \in \Xi$ be fixed. Using the Markov inequality, we obtain
\begin{eqnarray}
\Pr_{\theta,\xi} ( \Psi_{\alpha,D} =0)
& = & \Pr_{\theta,\xi} ( T_D < t_{1-\alpha,D}), \nonumber \\
& = & \Pr_{\theta,\xi} \left( T_D - \E_\theta[T_D] < t_{1-\alpha,D} - \sum_{k=1}^D \theta_k^2 \right) \nonumber \\
& \leq & \Pr_{\theta,\xi} \left( \left|T_D - \E_\theta[T_D]  \right| \geq  \sum_{k=1}^D \theta_k^2 - t_{1-\alpha,D}  \right)\nonumber \\
& \leq & \frac{ \mathrm{Var}_{\theta,\xi}(T_D)}{\left( \sum_{k=1}^D \theta_k^2 - t_{1-\alpha,D}  \right)^2},
\label{eq:2kind}
\end{eqnarray}
where we have implicitly assumed that 
$$ \sum_{k=1}^D \theta_k^2 > t_{1-\alpha,D}.$$
Now, we need an upper bound for the variance term. First remark that 
\begin{eqnarray}
\mathrm{Var}_\theta(T_D) 
& = & \mathrm{Var}_{\theta,\xi} \left( \sum_{k=1}^D b_k^{-2} (y_k^2 - \varepsilon^2)    \right)  \nonumber  \\
& = & \underbrace{ \sum_{k=1}^D b_k^{-4}  \mathrm{Var}_{\theta,\xi} ( y_k^2 -\varepsilon^2)}_{:= R_\theta(D)} +  \underbrace{\sum_{\substack{k,l=1 \\  k \not = l}}^D b_k^{-2} b_l^{-2} \mathrm{Cov}_{\theta,\xi}  \left( y_k^2-\varepsilon^2, y_l^2-\varepsilon^2 \right) }_{:=S_\theta(D)}.
\label{eq:var}
\end{eqnarray}
\underline{Calculation of $R_\theta(D)$:} Using simple algebra, we get, for all $k\in \mathcal{N}$,
\begin{eqnarray*}
\mathrm{Var}_{\theta}( y_k^2 -\varepsilon^2)
& = & \mathrm{Var}_{\theta,\xi}  \left[  ( b_k \theta_k + \varepsilon)^2 -\varepsilon^2  \right] \\
& = & \mathrm{Var}_{\theta,\xi}  \left[ b_k^2 \theta_k^2 + \varepsilon^2 \xi_k^2  + 2 b_k \theta_k \varepsilon \xi_k - \varepsilon^2 \right] \\
& = & \mathrm{Var}_{\theta,\xi}  \left[ \varepsilon^2 (\xi_k^2 -1)   + 2 b_k \theta_k \varepsilon \xi_k  \right] \\
& = & \varepsilon^4 \E [ (\xi_k^2-1)^2] + 4 \varepsilon^2 b_k^2 \theta_k^2 + 4 \varepsilon^3 b_k \theta_k  \E[\xi_k^3].
\end{eqnarray*}
Hence, using the last equality, we obtain 
\begin{eqnarray}
R_\theta(D) 
& \leq & C_1 \varepsilon^4 \sum_{k=1}^D b_k^{-4} + 4 \varepsilon^2 \sum_{k=1}^D b_k^{-2} \theta_k^2 + 4 C_2 \varepsilon^3 \sum_{k=1}^D b_k^{-3} | \theta_k |  \nonumber\\
& \leq & C_1 \varepsilon^4 \sum_{k=1}^D b_k^{-4}  + 4 \varepsilon^2 (\max_{1\leq k \leq D} b_k^{-2}) \sum_{k=1}^D \theta_k^2+   4 C_2 \varepsilon^3 \sum_{k=1}^D b_k^{-3} | \theta_k | ,
\label{eq:inter1} 
\end{eqnarray}
where the constant $C_1$ has been introduced in (\ref{eq:C_one}) and 
$$ C_2 := \sup_{\xi \in \Xi } \ \sup_{k\in \mathcal{N}} |\E[\xi_k^3]| < +\infty. $$
Note that, using first the Cauchy-Schwartz inequality and then the Peter-Paul inequality (see, e.g., \cite{CR_2016}, p. 18), we get 
\begin{eqnarray}
\varepsilon^3 \sum_{k=1}^D b_k^{-3} | \theta_k |
& = &  \sum_{k=1}^D \varepsilon^2  b_k^{-2}  \varepsilon b_k^{-1} | \theta_k | \nonumber  \\
& \leq & \sqrt{ \varepsilon^4 \sum_{k=1}^D   b_k^{-4}}  \sqrt{ \varepsilon^2 \sum_{k=1}^D b_k^{-2}  \theta_k^2} \nonumber \\
& \leq & \frac{1}{2} \varepsilon^4 \sum_{k=1}^D   b_k^{-4} + \frac{1}{2} \varepsilon^2 \sum_{k=1}^D b_k^{-2}  \theta_k^2  \nonumber \\
& \leq & \frac{1}{2} \varepsilon^4 \sum_{k=1}^D   b_k^{-4} + \frac{1}{2} \varepsilon^2 (\max_{1\leq k \leq D} b_k^{-2} )\sum_{k=1}^D  \theta_k^2.  
\label{eq:inter2}
\end{eqnarray}
Combining inequalities (\ref{eq:inter1}) and (\ref{eq:inter2}), we obtain 
\begin{equation}
R_\theta(D) \leq (C_1+ 2 C_2) \varepsilon^4 \sum_{k=1}^D b_k^{-4}  + (4+ 2 C_2) \varepsilon^2 (\max_{1\leq k \leq D} b_k^{-2}) \sum_{k=1}^D \theta_k^2.
\label{eq:H}
\end{equation}

\noindent
\underline{Calculation of $S_\theta(D)$:}  First, remark that, for all $k\in \mathcal{N}$, on noting that $ \E_{\theta,\xi}  [ (y_k^2 - \varepsilon^2)] = b_k^2 \theta_k^2$ and $\E[\xi_k^2]=1$, we get
\begin{eqnarray*}
\lefteqn{\mathrm{Cov}_{\theta,\xi}  \left( y_k^2-\varepsilon^2, y_l^2-\varepsilon^2 \right)}\\
& = & \E_{\theta,\xi}  \left[  (y_k^2 - \varepsilon^2 - b_k^2 \theta_k^2) \  (y_l^2 - \varepsilon^2 - b_l^2 \theta_l^2)   \right] \\
& = & \E_{\theta,\xi}  \left[  ((b_k\theta_k + \varepsilon \xi_k)^2 - \varepsilon^2 - b_k^2 \theta_k^2) \  ((b_l\theta_l + \varepsilon \xi_l)^2 - \varepsilon^2 - b_l^2 \theta_l^2)   \right]\\
& = & \varepsilon^4 \mathbb{E} [ (\xi_k^2-1) (\xi_l^2 -1)] + 4 \varepsilon^2 b_k b_l \theta_k \theta_l \mathbb{E} [ \xi_k \xi_l ]  + 2 \varepsilon^3 b_l \theta_l \E [ (\xi_k^2-1)\xi_l] + 2 \varepsilon^3 b_k \theta_k \E [ (\xi_l^2-1)\xi_k]. \\
\end{eqnarray*}
Hence, 
\begin{eqnarray}
S_\theta(D)
& = & 2\varepsilon^4 \sum_{\substack{k,l=1 \\  k \not = l}}^D  b_k^{-2} b_l^{-2}  \mathbb{E} [ (\xi_k^2-1) (\xi_l^2 -1)] \nonumber   \\
&  & + 8 \varepsilon^2 \sum_{\substack{k,l=1 \\  k \not = l}}^D b_k^{-1} b_l^{-1} \theta_k \theta_l \mathbb{E} [ \xi_k \xi_l ] \nonumber \\
& &  + 4  \varepsilon^3 \sum_{\substack{k,l=1 \\  k \not = l}}^D b_k^{-2} b_l^{-1}  \theta_l \E [ (\xi_k^2-1)\xi_l] \nonumber \\
& & + 4  \varepsilon^3 \sum_{\substack{k,l=1 \\  k \not = l}}^D b_l^{-2} b_k^{-1}  \theta_k \E [ (\xi_l^2-1)\xi_k].
\label{eq:S1}
\end{eqnarray}
Using the Cauchy-Schwarz inequality in each expectation of the above expression, we obtain 
\begin{eqnarray*}
S_\theta(D)
& \leq & 2 C_1 \varepsilon^4 \left(\sum_{k=1}^D b_k^{-2} \right)^2 + 8 \varepsilon^2 \left( \sum_{k=1}^D b_k^{-1} | \theta_k| \right)^2 + 8 \epsilon^3 C_1^{1/2} \left(  \sum_{k=1}^D b_k^{-1} |\theta_k| \right) \left( \sum_{k=1}^D b_k^{-2} \right).
\end{eqnarray*}
Then, using first the Peter-Paul inequality and then the Cauchy-Schwarz inequality, we get
\begin{eqnarray}
S_\theta(D) 
& \leq & 2( C_1 + 2C_1^{1/2})  \varepsilon^4 \left(\sum_{k=1}^D b_k^{-2} \right)^2 + 4(2+C_1^{1/2})  \varepsilon^2 \left( \sum_{k=1}^D b_k^{-1} | \theta_k| \right)^2  \nonumber \\
& \leq & 2( C_1 + 2C_1^{1/2})  \varepsilon^4 \left(\sum_{k=1}^D b_k^{-2} \right)^2 + 4(2+C_1^{1/2}) \varepsilon^2  \sum_{k=1}^D  b_k^{-2}   \sum_{k=1}^D \theta_k^2.
\label{eq:Psi}
\end{eqnarray}
Hence, combining (\ref{eq:var}), (\ref{eq:H}) and (\ref{eq:Psi}), we obtain 
\begin{eqnarray*} 
\mathrm{Var}_{\theta,\xi} (T_D)  
& \leq &  (C_1+ 2 C_2) \varepsilon^4 \sum_{k=1}^D b_k^{-4} \\
& & + 4( C_1 + 2C_1^{1/2})  \varepsilon^4 \left(\sum_{k=1}^D b_k^{-2} \right)^2 \\ 
& & + 2(10+  C_2+4C_1^{1/2}) \varepsilon^2 \sum_{k=1}^D  b_k^{-2} \sum_{k=1}^D \theta_k^2.
\end{eqnarray*}
For all $\gamma \in ]0,1[$, using again the Peter-Paul inequality, we get
\begin{eqnarray*}
\mathrm{Var}_{\theta,\xi} (T_D)  
& \leq &  (C_1+ 2 C_2) \varepsilon^4 \sum_{k=1}^D b_k^{-4} \\
& & + 4( C_1 + 2C_1^{1/2})  \varepsilon^4 \left(\sum_{k=1}^D b_k^{-2} \right)^2 \\ 
& & + \gamma^{-1}(10+  C_2+4C_1^{1/2}) \varepsilon^4 \left( \sum_{k=1}^D b_k^{-2}  \right)^2+ \gamma(10+  C_2+4C_1^{1/2})  \left( \sum_{k=1}^D \theta_k^2 \right)^2.
\end{eqnarray*}
Hence, since $\gamma \in ]0,1[$, it is easily seen that
\begin{eqnarray} 
\lefteqn{\mathrm{Var}_{\theta,\xi} (T_D) } \nonumber \\
& \leq & \gamma(10+  C_2+4C_1^{1/2})  \left( \sum_{k=1}^D \theta_k^2 \right)^2   \nonumber \\
& & +  \left\lbrace   (C_1+ 2 C_2) + 4( C_1 + 2C_1^{1/2}) + \gamma^{-1}(10+  C_2+4C_1^{1/2})   \right\rbrace\varepsilon^4 \left( \sum_{k=1}^D b_k^{-2}  \right)^2 \nonumber \\
& \leq & \gamma (10+  C_2+4C_1^{1/2})  \left( \sum_{k=1}^D \theta_k^2 \right)^2  +  \gamma^{-1}\left\lbrace   10 + 5C_1 + 2 C_2 +12C_1^{1/2}    \right\rbrace\varepsilon^4 \left( \sum_{k=1}^D b_k^{-2}  \right)^2 \nonumber \\
& \leq & ( 10 + 5C_1 + 2 C_2 +12C_1^{1/2}  ) \left\lbrace \gamma \left( \sum_{k=1}^D \theta_k^2 \right)^2  + \gamma^{-1}\varepsilon^4 \left( \sum_{k=1}^D b_k^{-2}  \right)^2\right\rbrace \nonumber \\
& := & K_2 \left\lbrace \gamma \left( \sum_{k=1}^D \theta_k^2 \right)^2  + \gamma^{-1}\varepsilon^4 \left( \sum_{k=1}^D b_k^{-2}  \right)^2 \right\rbrace.
\label{eq:var2}
\end{eqnarray}
Now, using (\ref{eq:2kind}) and (\ref{eq:var2}), and choosing $\gamma = \mathcal{C}_\beta^{-1}$, we get
\begin{eqnarray}
\Pr_{\theta,\xi}(\Psi_{\alpha,D} =0)
& \leq & \frac{ K_2 \left\lbrace \gamma \left( \sum_{k=1}^D \theta_k^2 \right)^2  + \gamma^{-1} \varepsilon^4 \left( \sum_{k=1}^D b_k^{-2}  \right)^2 \right\rbrace }{\left( \sum_{k=1}^D \theta_k^2 - K_1 \varepsilon^2  \sum_{k=1}^D b_k^{-2}    \right)^2} \nonumber \\
& \leq & \frac{ K_2 \left\lbrace \gamma + \gamma^{-1} \mathcal{C}_\beta^{-2} \right\rbrace }{\left( 1- K_1 \mathcal{C}_\beta^{-1}  \right)^2} \nonumber \\
& \leq & \frac{2 K_2 \mathcal{C}_\beta^{-1} }{\left( 1- K_1 \mathcal{C}_\beta^{-1}  \right)^2} \label{eq:Cb}\\
& \leq & \beta, \nonumber
\end{eqnarray}
provided that 
\begin{equation}
\sum_{k=1}^D  \theta_k^2 \geq  \mathcal{C}_\beta \varepsilon^2 \sum_{k=1}^D b_k^{-2} 
\label{eq:cond1}
\end{equation}
and $\mathcal{C}_\beta$ is the solution of the equation 
\begin{equation}
\frac{2 K_2 \mathcal{C}_\beta^{-1} }{\left( 1- K_1 \mathcal{C}_\beta^{-1}  \right)^2} = \beta.
\label{eq:condCb}
\end{equation}
To conclude the proof, since $\sum_{k> D} \theta_k^2 \leq a_D^{-2}$, remark that inequality (\ref{eq:cond1}) is satisfied provided that
$$ \| \theta \|^2 \geq \mathcal{C}_\beta \varepsilon^2 \sum_{k=1}^D b_k^{-2} + a_D^{-2} .$$
\begin{flushright}
$\Box$
\end{flushright}

\subsection{Proof of Proposition \ref{prop:lowerbound}}
\label{s:plowerbound}

When $\xi$ is Gaussian, we will write $\Xi_G$ instead of $\Xi$ and $\Pr_{\theta,\Sigma}$ instead of $\Pr_{\theta,\xi}$, where $\Sigma = (\Sigma_{kl})_{k,l \in \mathcal{N}}$ denotes the associated covariance matrix. We also define $S = \left\lbrace  \Sigma: \ \Sigma_{kk} =1 \right\rbrace$. Below, $\Psi_\alpha$ refers to an $\alpha$-level test. \\

Let $\theta^\star \in \Theta_a(r_\epsilon)$, $\xi^\star \in \Xi_G$ and $\Sigma^\star \in S$ be fixed. Their values will be made precise later on. Then
\begin{eqnarray*}
\inf_{\Psi_\alpha} \sup_{\substack{\theta \in \Theta_a(r_\epsilon) \\ \xi \in \Xi}} \Pr_{ \theta,\xi}(\Psi_\alpha =0)
& \geq & \inf_{\Psi_\alpha} \sup_{\substack{ \theta \in \Theta_a(r_\epsilon) \\ \xi \in \Xi_G}} \Pr_{ \theta,\xi}(\Psi_\alpha =0) \\
& = &  \inf_{\Psi_\alpha} \sup_{\substack{ \theta \in \Theta_a(r_\epsilon) \\ \tilde\Sigma \in S}} \Pr_{ \theta, \tilde \Sigma}(\Psi_\alpha =0) \\
& \geq &  \inf_{\Psi_\alpha} \sup_{\theta \in \Theta_a(r_\epsilon)} \Pr_{ \theta, \Sigma^\star}(\Psi_\alpha =0) \\
& \geq & \inf_{\Psi_\alpha}  \Pr_{ \theta^\star, \Sigma^\star}(\Psi_\alpha =0) \\
& \geq & 1-\alpha - \frac{1}{2}\left(\E_0[ L^2_{\theta^\star,\Sigma^\star}(Y)]  -1 \right)^{1/2},
\end{eqnarray*}
where $ L_{\theta^\star,\Sigma^\star}(Y) = d\Pr_{\theta^\star,\Sigma^\star}(Y)/d\Pr_{0,\Sigma^\star}(Y) $ is the likelihood ratio between the probability measures $\Pr_{\theta^\star,\Sigma^\star}$ and $\Pr_{0,\Sigma^\star}$ (for the last inequality, we refer to, e.g., (3.1) in \cite{MS_2014}). In particular, if we can find $\theta^\star$ and $\Sigma^\star$ such that 
$$
\E_0[ L^2_{\theta^\star,\Sigma^\star}(Y)] \leq  \mathcal{C}_{\alpha,\beta}
$$
for some $\beta \in ]0, 1-\alpha[$, then 
$$ \inf_{\Psi_\alpha} \sup_{\substack{\theta \in \Theta_a(r_\epsilon) \\ \xi \in \Xi}} \Pr_{ \theta,\xi}(\Psi_\alpha =0) \geq \beta.$$

\vspace{0.5cm}

Let $D \in \mathcal{N}$ be fixed. Now, we impose the following conditions on $\theta^\star$ and $\Sigma^\star$:
$$ \theta_k^\star =0 \quad \forall k > D \quad \mathrm{and} \quad \Sigma^\star_{kl} = 0 \quad \forall k>D,\  l>D, \ k\not = l.$$
Let $\Sigma_D^\star =(\Sigma^\star_{kl})_{1\leq k,l \leq D}$ be the remaining submatrix of $\Sigma^\star$. With a slight abuse of notation, we denote below $Y=(Y_1,\dots, Y_D)'$ and $b\theta^\star =(b_1 \theta^\star_1,\dots, b_D \theta^\star_D)'$. Then, by simple algebra,
\begin{eqnarray*}
L_{\theta^\star,\Sigma^\star}(Y) 
& = & \frac{\exp\left( -\frac{1}{2\varepsilon^2} (Y-b\theta^\star)' (\Sigma_D^\star)^{-1} (Y- b\theta^\star)    \right)}{\exp\left( -\frac{1}{2\varepsilon^2} Y' (\Sigma_D^\star)^{-1} Y \right)  } \\
& = & \exp \left[  \frac{1}{2\varepsilon^2}  \left\lbrace 2 (b\theta^\star)' (\Sigma_D^\star)^{-1} Y - (b\theta^\star)' (\Sigma_D^\star)^{-1} b\theta^\star    \right\rbrace \right].
\end{eqnarray*}
Hence, 
\begin{eqnarray*}
\E_0[ L^2_{\theta^\star,\Sigma^\star}(Y)] 
& = & \exp \left(-\frac{1}{\varepsilon^2} (b\theta^\star)' (\Sigma_D^\star)^{-1} b\theta^\star \right) \E_0 \left[  \exp\left( \frac{2}{\varepsilon^2} (b\theta^\star)' (\Sigma_D^\star)^{-1} Y \right)\right].
\end{eqnarray*}
It is easily seen that
\begin{eqnarray*}
\lefteqn{\E_0 \left[  \exp\left( \frac{2}{\varepsilon^2} (b\theta^\star)' (\Sigma_D^\star)^{-1} Y \right)\right] } \\
& = & \frac{1}{(2\pi \varepsilon^2)^{D/2} | \Sigma_D^\star |^{1/2}}  \int_{\mathbb{R}^D}   \exp\left( \frac{2}{\varepsilon^2} (b\theta^\star)' (\Sigma_D^\star)^{-1} y \right) \exp \left(  -\frac{1}{2\epsilon^2} y'( \Sigma_D^\star)^{-1} y \right) dy \\
& = & \frac{1}{(2\pi \varepsilon^2)^{D/2} | \Sigma_D^\star |^{1/2}}  \int_{\mathbb{R}^D}  \exp \left(  -\frac{1}{2\epsilon^2} \left\lbrace   y'( \Sigma_D^\star)^{-1} y  -4(b\theta^\star)' (\Sigma_D^\star)^{-1} y  \right\rbrace\right) dy \\
& = & \exp \left( \frac{2}{\varepsilon^2} (b\theta^\star)' (\Sigma_D^\star)^{-1} b\theta^\star  \right) \frac{1}{(2\pi \varepsilon^2)^{D/2} | \Sigma_D^\star |^{1/2}}  \int_{\mathbb{R}^D}  \exp \left(  -\frac{1}{2\epsilon^2}   (y-2b\theta^\star)'( \Sigma_D^\star)^{-1} (y-2b\theta^\star) \right) dy \\
& = & \exp \left( \frac{2}{\varepsilon^2} (b\theta^\star)' (\Sigma_D^\star)^{-1} b\theta^\star  \right).
\end{eqnarray*}
Hence, 
$$ \E_0[ L^2_{\theta^\star,\Sigma^\star}(Y)]  =\exp \left( \frac{1}{\varepsilon^2} (b\theta^\star)' (\Sigma_D^\star)^{-1} b\theta^\star  \right).$$
Now, we select $ \theta^\star$ as follows 
$$ \theta^\star_k = \frac{r_\varepsilon b_k^{-1} (\Sigma_D^\star v)_k}{\rho} \ \ \forall k\in \lbrace 1,\dots, D \rbrace \quad \mathrm{and} \quad \theta_k^\star = 0 \ \ \forall k> D,$$
where 
\begin{equation}
\rho^2 = \sum_{k=1}^D b_k^{-2} (\Sigma_D^\star v)_k^2 \quad  \mathrm{and} \quad v_k = \frac{1}{\sqrt{D}} \; \forall k\in \lbrace 1,\dots, D \rbrace. 
\label{eq:rho2}
\end{equation}
Now, define $v =(v_1,\dots, v_D)'$ and note that $\| v \| =1$. Then, it is easily seen that 
\begin{equation}
\E_0[ L^2_{\theta^\star,\Sigma^\star}(Y)]  = \exp \left( \frac{r_\varepsilon^2}{\varepsilon^2 \rho^2} v'   \Sigma_D^\star v  \right).
\label{eq:E}
\end{equation}

We first construct a specific $\xi^\star \in \Xi_G$. 
Let 
\begin{equation}
\xi_k^\star = d_k \, \eta_0 + \sqrt{1- d_k^2} \, \eta_k \ \  \forall k \in \lbrace 1,\dots, D \rbrace \quad \mathrm{and} \quad \xi_k^\star = \eta_k \ \ \forall k> D,
\label{eq:structure}
\end{equation}
where $(\eta_k)_{k\in \mathcal{N}_0}$ denotes a sequence of independent standard Gaussian random variables and $d = (d_k)_{1\leq k \leq D}$ is a real sequence such that $1/\sqrt{2} \leq d_k <1 $ for all $1\leq k \leq D$. Obviously, $\xi^\star \in \Xi_G$, since $\xi^\star$ is Gaussian,
$$ \E[\xi_k^\star] =0, \; \E[(\xi_k^\star)^2] =1 \quad \forall k \in \mathcal{N}$$ 
and
$$  \max_{1\leq k \leq D} \E[(\xi_k^\star)^4] \leq \max_{1 \leq k \leq D} \left[8 d_k^4 \E[ \xi_0^4] + 8 (1-d_k)^2 \E[\eta_k^4]  \right] \leq 16 \E [\xi_0^4] =C < +\infty,$$
$$ \sup_{k> D} \E[(\xi_k^\star)^4] = 3.$$
Now, we need to bound the expression in (\ref{eq:E}). Using (\ref{eq:structure}), we get
$$ \Sigma_{kl}^\star = \E [ \xi_k^\star \xi_l^\star ] = d_k d_l \geq 1/2 \ \ \forall k,l \in \lbrace 1,\dots, D \rbrace \; \mathrm{with} \; k\not = l \quad \mathrm{and} \quad \Sigma_{kl} = 0 \ \ \forall k,l > D, \ k\not = l.$$
Note also that, since $v$ is a unit vector,
$$ v' \Sigma_D^\star v  \leq \max_{\| a \| =1} a' \Sigma_D^\star a  \leq D$$
since the largest eigenvalue of $\Sigma_D^\star$ is smaller than $D$.  Now, using (\ref{eq:rho2}), we get 
$$ \rho^2 = \sum_{k=1}^D b_k^{-2} (\Sigma_D^\star v)_k^2 = \sum_{k=1}^D \left(  \sum_{l=1}^D \Sigma^\star_{kl} v_l \right)^2 \geq \frac{1}{4} D \sum_{k=1}^D b_k^{-2}$$
since $\Sigma^\star_{kl} \geq 1/2$ for all $k,l \in \lbrace 1,\dots, D \rbrace$. Hence,
$$  \E_0[ L^2_{\theta^\star,\Sigma^\star}(Y)] \leq \exp \left( \frac{4r_\varepsilon^2}{ \varepsilon^2} \frac{1}{\sum_{k=1}^D b_k^{-2}}   \right)  \leq \mathcal{C}_{\alpha,\beta}$$
provided
$$  \| \theta \|^2= r_\varepsilon^2 \leq  \left( \frac{1}{4} \ln( \mathcal{C}_{\alpha,\beta}) \right)\varepsilon^2 \sum_{k=1}^D b_k^{-2}.$$ 
To conclude the proof, we need to ensure that the constructed $\theta^\star$ belongs to $\mathcal{E}_a$. Remark that, since $a$ is an increasing sequence, 
\begin{eqnarray*}
\sum_{k \in \mathcal{N}} a_k^2 (\theta_k^\star)^2 = \sum_{k=1}^D a_k^2 (\theta_k^\star)^2 
 \leq  a_D^2 \sum_{k=1}^D a_k^2 (\theta_k^\star)^2 
 \leq  a_D^2 r_\varepsilon^2 
 \leq  1
\end{eqnarray*}
provided $r_\varepsilon^2 \leq a_D^{-2}$. Hence, 
$$ \inf_{\Psi_\alpha} \sup_{\substack{\theta \in \Theta_a(r_{\varepsilon}) \\ \xi \in \Xi}} \Pr_{ \theta,\xi}(\Psi_\alpha =0) \geq \beta,$$
as soon as 
$$ r_\varepsilon^2 \leq \left( \frac{1}{4} \ln( \mathcal{C}_{\alpha,\beta}) \right)\varepsilon^2 \sum_{k=1}^D b_k^{-2} \wedge a_D^{-2}.$$
\begin{flushright}
$\Box$
\end{flushright}

\subsection{Proof of Theorem \ref{thm:1}}
\label{s:pthm1}

In Proposition \ref{prop:upperbound}, we have proved that for all $D \in \mathcal{N}$, there exists an $\alpha$-level test $\Psi_{\alpha,D}$ such that
$$ \sup_{\substack{\theta\in \Theta_a(r_{\varepsilon,D}) \\ \xi \in \Xi}} \Pr_{\theta,\xi}(\Psi_{\alpha,D}=0) \leq \beta,$$
for all radius $r_{\varepsilon,D}>0$ satisfying
$$ r^2_{\varepsilon,D} \geq \mathcal{C}_\beta \varepsilon^2 \sum_{k=1}^D b_k^{-2} + a_D^{-2}.$$ 
Now, setting
$$ D^\dagger = \arg \inf_{D\in \mathcal{N}} \left[ \mathcal{C}_\beta   \varepsilon^2 \sum_{k=1}^D b_k^{-2} + a_D^{-2}\right]$$
and denoting by $ \Psi_{\alpha,D^\dagger}$ the associated $\alpha$-level test, we get
$$ \sup_{\substack{\theta\in \Theta_a(r_\varepsilon^\star) \\ \xi \in \Xi}} \Pr_{\theta,\xi}(\Psi_{\alpha,D^\dagger}=0) \leq \beta,$$
for all radius $r_{\varepsilon}^\star>0$ satisfying
$$ (r_{\varepsilon}^\star)^2 \geq \inf_{D\in \mathcal{N}} \left[\mathcal{C}_\beta \varepsilon^2 \sum_{k=1}^D b_k^{-2} + a_D^{-2}\right].$$ 
Hence,
$$ \tilde r_\varepsilon^2 \leq \inf_{D\in\mathcal{N}} \left[ \mathcal{C}_\beta   \varepsilon^2 \sum_{k=1}^D b_k^{-2} + a_D^{-2}\right] .$$
Similarly, using Proposition \ref{prop:lowerbound}, 
$$ \inf_{\Psi_\alpha} \sup_{\substack{\theta \in \Theta_a(r_{\varepsilon}) \\ \xi \in \Xi}} \Pr_{ \theta,\xi}(\Psi_\alpha =0) \geq \beta,$$
for all radius $r_{\varepsilon}>0$ such that 
$$ r_\varepsilon^2 \leq \left( \frac{1}{4} \ln( \mathcal{C}_{\alpha,\beta}) \right)\varepsilon^2 \sum_{k=1}^D b_k^{-2} \wedge a_D^{-2}.$$
This results occurs for all $D \in \mathcal{N}$. Hence, 
$$ \inf_{\Psi_\alpha} \sup_{\substack{\theta \in \Theta_a(r_{\varepsilon}) \\ \xi \in \Xi}} \Pr_{ \theta,\xi}(\Psi_\alpha =0) \geq \beta,$$
for all radius $r_{\varepsilon}>0$ such that 
$$ r_\varepsilon^2 \leq \sup_{D\in \mathcal{N}} \left[ \left( \frac{1}{4} \ln( \mathcal{C}_{\alpha,\beta}) \right)\varepsilon^2 \sum_{k=1}^D b_k^{-2} \wedge a_D^{-2}\right].$$
This entails that 
$$ \tilde r_\varepsilon^2 \geq \sup_{D\in \mathcal{N}} \left[ \left( \frac{1}{4} \ln( \mathcal{C}_{\alpha,\beta}) \right)\varepsilon^2 \sum_{k=1}^D b_k^{-2} \wedge a_D^{-2}\right].$$

\begin{flushright}
$\Box$
\end{flushright}

\subsection{Proof of Proposition \ref{prop:short1}}
\label{s:pshort1}

Remark that, under Assumption $\mathcal{H}_\mathcal{D}$, 
\begin{eqnarray*}
\frac{S_0(D)}{\varepsilon^4}
& \sim & \sum_{\substack{k,l=1 \\  k \not = l}}^D k^{2t} l^{2t} | k-l |^{-2s} \\
& = & 2\sum_{k=1}^D k^{2t} \sum_{l=1}^{k-1} \frac{l^{2t}}{(k-l)^{2s}}\\
& \leq & 2 \sum_{k=1}^D k^{4t} \sum_{l=1}^{k-1} \frac{1}{(k-l)^{2s}}. 
\end{eqnarray*}
Then, for all $k\in \mathcal{N}$, 
$$ \sum_{l=1}^{k-1} \frac{1}{(k-l)^{2s}}  = \sum_{m=1}^{k-1} \frac{1}{m^{2s}}  \sim k^{1-2s}.$$
In particular, the above sum is finite whatever the value of $k$ provided $s > 1/2$. Hence, under Assumption $\mathcal{H}_\mathcal{D}$, we get
$$ S_0(D) \lesssim \varepsilon^4 \sum_{k=1}^D k^{4t} k^{1-2s} \sim \varepsilon^4 D^{4t - 2s +2} = o(R_0(D)) \quad \mathrm{as} \ D\rightarrow +\infty.$$

\begin{flushright}
$\Box$
\end{flushright}

\subsection{Proof of Proposition \ref{prop:short2}}
\label{s:pshort2}

Recall from (\ref{eq:var}) that 
$$
\mathrm{Var}_{\theta,\xi} (T_D)  =  R_\theta(D) +  S_\theta(D), 
 $$
 where, using (\ref{eq:H}),
\begin{eqnarray*}
R_\theta(D) & = & \sum_{k=1}^D b_k^{-4}  \mathrm{Var}_{\theta,\xi} ( y_k^2 -\varepsilon^2) \\
& \lesssim & \varepsilon^4 \sum_{k=1}^D b_k^{-4}  + \varepsilon^2 (\max_{1\leq k \leq D} b_k^{-2}) \sum_{k=1}^D \theta_k^2.
\end{eqnarray*}
Moreover, using (\ref{eq:S1}),
\begin{eqnarray*}
S_\theta(D)
& =& \sum_{\substack{k,l=1 \\  k \not = l}}^D b_k^{-2} b_l^{-2} \mathrm{Cov}_{\theta,\xi}  \left( y_k^2-\varepsilon^2, y_l^2-\varepsilon^2 \right)\\
& = & 2\varepsilon^4 \sum_{\substack{k,l=1 \\  k \not = l}}^D b_k^{-2} b_l^{-2}  \mathbb{E} [ (\xi_k^2-1) (\xi_l^2 -1)]  \\
&  & + 8 \varepsilon^2 \sum_{\substack{k,l=1 \\  k \not = l}}^D b_k^{-1} b_l^{-1} \theta_k \theta_l \mathbb{E} [ \xi_k \xi_l ]  \\
& &  + 4  \varepsilon^3 \sum_{\substack{k,l=1 \\  k \not = l}}^D b_k^{-2} b_l^{-1}  \theta_l \E [ (\xi_k^2-1)\xi_l]  \\
& & + 4  \varepsilon^3 \sum_{\substack{k,l=1 \\  k \not = l}}^D b_l^{-2} b_k^{-1}  \theta_k \E [ (\xi_l^2-1)\xi_k] \\
& := & R_1 + R_2 + R_3 + R_4.
\end{eqnarray*}
Note that, using the above proposition, $R_1 = S_0(D) = o(R_\theta(D))$ as $D \rightarrow + \infty$. Then, using (\ref{eq:anoise2}), we can immediately see that $R_3 = R_4 = 0$. In order to conclude, using the Cauchy-Schwarz  and Peter-Paul inequalities, we get, for any $\gamma \in ]0,1[$,
\begin{eqnarray*}
R_2 & = & 8 \varepsilon^2 \sum_{\substack{k,l=1 \\  k \not = l}}^D b_k^{-1} b_l^{-1} \theta_k \theta_l \mathrm{Cov}(\xi_k, \xi_l ) \\
& \lesssim &   \varepsilon^2 \sum_{\substack{k,l=1 \\  k \not = l}}^D b_k^{-1} b_l^{-1} |\theta_k|   |\theta_l|  \frac{1}{| k-l |^s}   \\
& \leq & \sqrt{\varepsilon^4 \sum_{\substack{k,l=1 \\  k \not = l}}^D b_k^{-2} b_l^{-2} \frac{1}{|k-l |^{2s}}} \sqrt{ \sum_{\substack{k,l=1 \\  k \not = l}}^D \theta_k^2   \theta_l^2 }   \\
& \leq & \gamma \left( \sum_{k=1}^D \theta_k^2 \right)^2 + \gamma^{-1} \varepsilon^4 \sum_{\substack{k,l=1 \\  k \not = l}}^D b_k^{-2} b_l^{-2} \frac{1}{| k-l |^{2s}} \\
& \sim & \gamma \left( \sum_{k=1}^D \theta_k^2 \right)^2 + \gamma^{-1} R_1. 
\end{eqnarray*}
Summarizing all the above computations, we obtain, for any $\gamma \in ]0,1[$,
\begin{eqnarray*}
\mathrm{Var}_{\theta,\xi} (T_D) 
& \lesssim & (1+\gamma^{-1}) \varepsilon^4 \sum_{k=1}^D b_k^{-4}  + \varepsilon^2 (\max_{1\leq k \leq D} b_k^{-2}) \sum_{k=1}^D \theta_k^2 + \gamma \left( \sum_{k=1}^D \theta_k^2 \right)^2 \\
& \lesssim & (1+\gamma^{-1}) \varepsilon^4 \sum_{k=1}^D b_k^{-4}  + \gamma \left( \sum_{k=1}^D \theta_k^2 \right)^2,
\end{eqnarray*}
where we have used again the Peter-Paul inequality (see, e.g., \cite{CR_2016}, p. 18).

\begin{flushright}
$\Box$
\end{flushright}

\bigskip
\bigskip
\bibliography{Survey}

\begin{thebibliography}{10}

\bibitem{Baraud}
Y.~Baraud.
\newblock Non-asymptotic minimax rates of testing in signal detection.
\newblock {\em Bernoulli}, 8(5):577--606, 2002.

\bibitem{Munk_testing}
N.~Bissantz, G.~Claeskens, H.~Holzmann, and A.~Munk.
\newblock Testing for lack of fit in inverse regression---with applications to
  biophotonic imaging.
\newblock {\em Journal of the Royal Statistical Society, Series B,},
  71(1):25--48, 2009.

\bibitem{CR_2016}
R.E. Castillo and H.~Rafeiro.
\newblock {\em An Introductory Course in {L}ebesgue Spaces}.
\newblock CMS Books in Mathematics/Ouvrages de Math\'ematiques de la SMC.
  Springer, 2016.

\bibitem{Cavalier_frac}
L.~Cavalier.
\newblock Estimation in a problem of fractional integration.
\newblock {\em Inverse Problems}, 20(5):1445--1454, 2004.

\bibitem{Cavalier_book}
L.~Cavalier.
\newblock Inverse problems in statistics.
\newblock In {\em Inverse Problems and High-Dimensional Estimation}, volume 203
  of {\em Lect. Notes Stat. Proc.}, pages 3--96. Springer, Heidelberg, 2011.

\bibitem{ISS_2011}
Yu.I. Ingster, T.~Sapatinas, and I.A. Suslina.
\newblock Minimax nonparametric testing in a problem related to the {R}adon
  transform.
\newblock {\em Mathematical Methods of Statistics}, 20(4):347--364, 2011.

\bibitem{ISS_2012}
Yu.I. Ingster, T.~Sapatinas, and I.A. Suslina.
\newblock Minimax signal detection in ill-posed inverse problems.
\newblock {\em Annals of Statistics}, 40:1524--1549, 2012.

\bibitem{IS_2003}
Yu.I. Ingster and I.A. Suslina.
\newblock {\em Nonparametric Goodness-of-Fit Testing Under Gaussian Models},
  volume 169 of {\em Lecture Notes in Statistics}.
\newblock Springer-Verlag, New York, 2003.

\bibitem{Isserlis}
L.~Isserlis.
\newblock On a formula for the product-moment coefficient of any order of a
  normal frequency distribution in any number of variables.
\newblock {\em Biometrika}, 12(1/2):134--139, 1918.

\bibitem{J_1999}
I.M. Johnstone.
\newblock Wavelet shrinkage for correlated data and inverse problems:
  adaptivity results.
\newblock {\em Statistica Sinica}, 9(1):51--83, 1999.

\bibitem{LLM_2011}
B.~Laurent, J.-M. Loubes, and C.~Marteau.
\newblock Testing inverse problems: a direct or an indirect problem?
\newblock {\em Journal of Statistical Planning and Inference},
  141(5):1849--1861, 2011.

\bibitem{LLM_2012}
B.~Laurent, J.-M. Loubes, and C.~Marteau.
\newblock Non asymptotic minimax rates of testing in signal detection with
  heterogeneous variances.
\newblock {\em Electronic Journal of Statistics}, 6:91--122, 2012.

\bibitem{Mallat}
S.~Mallat.
\newblock {\em A Wavelet Tour of Signal Processing. 2nd Edition.}
\newblock Academic Press, San Diego, 1999.

\bibitem{MM_2013}
C.~Marteau and P.~Math{\'e}.
\newblock General regularization schemes for signal detection in inverse
  problems.
\newblock {\em Mathematical Methods of Statistics}, 23(3):176--200, 2014.

\bibitem{MS_2014}
C.~Marteau and T.~Sapatinas.
\newblock A unified treatment for non-asymptotic and asymptotic approaches to
  minimax signal detection.
\newblock {\em Statistics Surveys}, 9:253--297, 2015.

\bibitem{Tsybakov}
A.B. Tsybakov.
\newblock {\em Introduction to Nonparametric Estimation}.
\newblock Springer Series in Statistics. Springer, New York, 2009.
\newblock Revised and extended from the 2004 French original, Translated by
  Vladimir Zaiats.

\end{thebibliography}
\bibliographystyle{plain}

\end{document}